\documentclass[a4paper,12pt]{article}
\usepackage[T1]{fontenc}
\usepackage[utf8]{inputenc}
\usepackage[english]{babel}

\usepackage{amsmath, latexsym, amsfonts, amssymb}
\usepackage{amsthm}
\usepackage{mathrsfs,enumerate}
\usepackage{graphics,epsfig,color}
\usepackage[top=1.4cm, bottom=1.7cm, left=2.1cm, right=2.8cm]{geometry}
\usepackage{enumitem}
\usepackage[svgnames]{xcolor}
\usepackage{comment}

\usepackage{fullpage}
\usepackage{hyperref}
\hypersetup{
    colorlinks=true,
    linkcolor=OrangeRed,
    citecolor=Goldenrod,
    urlcolor=DarkGreen,
    pdfborder={0 0 0}
}
\definecolor{fond}{rgb}{0.05,0.05,0.25}

\DeclareMathOperator{\Fix}{Fix}
\DeclareMathOperator{\Poiss}{Poiss}
\DeclareMathOperator{\dvt}{d_{TV}}

\title{Limit profile for random transpositions}
\author{Lucas Teyssier\footnote{Pronunciation: [lyka tesje]. Student at ENS and Sorbonne Université. Current email adress: \texttt{lucas.teyssier@ens.fr}}}
\date{May 2019}

\begin{document}

\maketitle

\theoremstyle{definition}
\newtheorem{thm}{Theorem}[section]
\newtheorem{problem}{Problem}[section]
\newtheorem{notation}[thm]{Notation}
\newtheorem{claim}[thm]{Claim}
\newtheorem{lemma}[thm]{Lemma}
\newtheorem{corollary}[thm]{Corollary}
\newtheorem{prop}[thm]{Proposition}
\newtheorem{defn}[thm]{Definition}
\newtheorem{conjecture}[thm]{Conjecture}
\newtheorem{remark}[thm]{Remark}

\newcommand{\cA}{{\ensuremath{\mathcal A}} }
\newcommand{\cB}{{\ensuremath{\mathcal B}} }
\newcommand{\cC}{{\ensuremath{\mathcal C}} }
\newcommand{\cD}{{\ensuremath{\mathcal D}} }
\newcommand{\cE}{{\ensuremath{\mathcal E}} }
\newcommand{\cF}{{\ensuremath{\mathcal F}} }
\newcommand{\cG}{{\ensuremath{\mathcal G}} }
\newcommand{\cH}{{\ensuremath{\mathcal H}} }

\newcommand{\cL}{{\ensuremath{\mathcal L}} }
\newcommand{\cM}{{\ensuremath{\mathcal M}} }
\newcommand{\cN}{{\ensuremath{\mathcal N}} }
\newcommand{\cO}{{\ensuremath{\mathcal O}} }
\newcommand{\cP}{{\ensuremath{\mathcal P}} }

\newcommand{\bbE}{{\ensuremath{\mathbb E}} }
\newcommand{\bbN}{{\ensuremath{\mathbb N}} }
\newcommand{\bbZ}{{\ensuremath{\mathbb Z}} }
\newcommand{\bbQ}{{\ensuremath{\mathbb Q}} }
\newcommand{\bbR}{{\ensuremath{\mathbb R}} }
\newcommand{\bbC}{{\ensuremath{\mathbb C}} }
\newcommand{\bbP}{{\ensuremath{\mathbb P}} }

\newcommand{\rstel}{\textcolor{red}{*}}
\newcommand{\rs}{\textcolor{red}{*}}
\newcommand{\nn}{{\ensuremath{n \in \mathbb N}}}
\newcommand{\kn}{{\ensuremath{k \in \mathbb N}}}
\newcommand{\bigcupnn}{{\ensuremath{\bigcup_{n \in \mathbb N}}}}
\newcommand{\bigcapnn}{{\ensuremath{\bigcap_{n \in \mathbb N}}}}

\newcommand{\limn}{{\ensuremath{\lim_{n \to \infty \,}}}}
\newcommand{\limcn}{{\ensuremath{\lim_{n \to \infty} \uparrow \,}}}
\newcommand{\limdn}{{\ensuremath{\lim_{n \to \infty} \downarrow \,}}}
\newcommand{\limck}{{\ensuremath{\lim_{k \to \infty} \uparrow \,}}}
\newcommand{\limdk}{{\ensuremath{\lim_{k \to \infty} \downarrow \,}}}
\newcommand{\sumnn}{{\ensuremath{\sum_{n \in \mathbb N}}}}
\newcommand{\inv}{{\ensuremath{^{-1}}}}

\newcommand{\ag}{\left\{ } 
\newcommand{\cog}{\left]} 
\newcommand{\cfg}{\left[}
\newcommand{\ad}{\right\} }
\newcommand{\cod}{\right[}
\newcommand{\cfd}{\right]}
\newcommand{\pg}{\left(} 
\newcommand{\pd}{\right)} 
\newcommand{\bbRb}{\ensuremath{\overline{\bbR}}} 
\newcommand{\npos}{{\ensuremath{n \geq 0}}} 
\newcommand{\du}{{\ensuremath{\;:\;}}} 
\newcommand{\rbr}{{\ensuremath{(\bbR,\cB(\bbR))}}}
\newcommand{\ea}{{\ensuremath{(E,\cA)}}}
\newcommand{\eamu}{{\ensuremath{(E,\cA, \mu)}}}
\newcommand{\fb}{{\ensuremath{(F,\cB)}}}
\newcommand{\ifm}{{\ensuremath{\int fd\mu}}}
\newcommand{\axe}{\ag x\in E\du}
\newcommand{\pch}{\text{ p.ĉ.}}
\newcommand{\pc}{\text{ p.c.}}
\newcommand{\la}{\ensuremath{\lambda}}
\newcommand{\ch}{\text{ch}}

\newcommand{\summ}[2]{\sum_{#1}^{#2}}
\newcommand{\prodd}[2]{\prod_{#1}^{#2}}

\newcommand{\intmu}[1]{\int #1 d\mu}
\newcommand{\intmup}[1]{\int \pg #1 \pd d\mu}
\newcommand{\lequ}[2]{{{#1} \leq {#2}}}
\newcommand{\gequ}[2]{{{#1} \geq {#2}}}
\newcommand{\abs}[1]{\left\lvert #1 \right\rvert}
\newcommand{\absa}[1]{\bigl\lvert #1 \bigr\rvert}
\newcommand{\absb}[1]{\Bigl\lvert #1 \Bigr\rvert}
\newcommand{\absc}[1]{\biggl\lvert #1 \biggr\rvert}
\newcommand{\absd}[1]{\Biggl\lvert #1 \Biggr\rvert}

\newcommand{\ura}[1]{\underset{#1}{\longrightarrow}}
\newcommand{\norme}[1]{|| #1 ||}

\newcommand{\Sn}{\ensuremath{\mathfrak{S}_n}}
\newcommand{\Pek}{\ensuremath{P^{*k}}}
\newcommand{\Pet}{\ensuremath{P^{*t}}}

\newcommand{\nl}{\ensuremath{\frac{1}{2}n\log(n)}}
\newcommand{\nlc}{\ensuremath{\frac{1}{2}n\log(n)+cn}}

\setcounter{tocdepth}{2} 
\begin{center}
    \textbf{Résumé\footnote{Cet article possède aussi une version en français.}}
\end{center}
Nous présentons une amélioration du lemme de majoration de Diaconis, qui permet de calculer la valeur limite de la distance à la stationnarité. Nous l'appliquons ensuite aux transpositions aléatoires étudiées par Diaconis et Shahshahani.
\begin{center}
    \textbf{Rezumo\footnote{Tiu artikolo ankaŭ havas version en Esperanto.}}
\end{center}
Ni prezentas plibonigon de la superbara lemo de Diaconis, kiu ebligas nin kalkuli la limesan valoron de la distanco al staranteco. Ni poste aplikas ĝin al hazardaj $2$-cikloj studitaj de Diaconis kaj Shahshahani.
\begin{center}
    \textbf{Abstract}
\end{center}
We present an improved version of Diaconis' upper bound lemma, which is used to compute the limiting value of the distance to stationarity. We then apply it to random transpositions studied by Diaconis and Shahshahani. 


\tableofcontents
\section{Introduction}\label{section1}

\subsection{Main results}\label{subsection11}

Let $\Sn$ be the symmetric group of indice $n$ and $P_n$ the probability on $\Sn$ defined by 
\begin{equation*}
    P_n(Id) = \frac{1}{n} \text{ and } P_n(\tau)=\frac{2}{n^2} \text{ if } \tau \text{ is a transposition.}
\end{equation*}
This is the random transposition shuffle on $\Sn$, as studied in a landmark paper of Diaconis and Shashahani \cite{DiaconisShahshahani1981}. \\
Let also $U_n$ be the uniform probability on $\Sn$.
If $E$ is a set and $\mu$, $\nu$ are probabilities on $E$, we define the total variation distance\footnote{In the proofs we will use the $L^1$ distance, noted $\text{d}_1$, in order not to carry the factor $\frac{1}{2}$.} between $\mu$ and $\nu$ by the formula
\begin{equation*}
    \dvt(\mu, \nu) = \frac{1}{2} \text{d}_{\text{1}}(\mu, \nu) = \frac{1}{2}\sum_{x\in E}\abs{\mu(x)-\nu(x)}.
\end{equation*}
In \cite{DiaconisShahshahani1981}, Diaconis and Shahshahani showed that this random walk undergoes a cutoff phenomenon at $\frac{1}{2}n\log(n)$, i.e., letting $f(n) = \left\lfloor\frac{1}{2}n\log(n)\right\rfloor$, that for all $0<\epsilon<1$,
\begin{equation*}
    \dvt\left(P_n^{*(1-\epsilon)f(n)},U_n\right) \xrightarrow[n\rightarrow\infty]{} 1 \;\;\;\;\;\text{ and }\;\;\;\;\; \dvt\left(P_n^{*(1+\epsilon)f(n)},U_n\right) \xrightarrow[n\rightarrow\infty]{} 0.
\end{equation*}
Despite a lot of work on mixing times in general and on random transpositions in particular (see references below), obtaining a precise description of the way this transition occurs has remained an open problem, formally asked by Nathanaël Berestycki at an AIM workshop on Markov chains mixing times in 2016 (\url{http://aimpl.org/markovmixing/5/}).  \\ 
Our main result is the following:
\begin{thm}\label{theorem}
Let $c\in\bbR$.  Then we have: 
\begin{equation*}
    \dvt\left(P_n^{*\left\lfloor\frac{1}{2}n\log(n)+cn\right\rfloor}, U_n\right) \xrightarrow[n\rightarrow\infty]{} \dvt\left(\Poiss\left(1+e^{-2c}\right),\Poiss(1)\right),
\end{equation*}
where $\text{Poiss}(a)$ stands for the Poisson law of parameter $a$.\\
\end{thm}
\paragraph{Limiting profile conjectures} \label{remarktheorem} \textbf{ } \\
We anticipate the limiting profile $\dvt\left(\text{Poiss}\left(1+e^{-c}\right),\text{Poiss}(1)\right)$, which we obtain in our problem if we replace the time $\left\lfloor\frac{1}{2}n\log(n)+cn\right\rfloor$ by a slightly more natural time, $\left\lfloor\frac{1}{2} (n\log(n)+cn)\right\rfloor$, to arise for many other mixing time problems on $\Sn$, namely the problems where the last things to be mixed are the fixed points. It seems to be often the case when the probability $P_n$ is constant over conjugacy classes. For example, using the formulas in \cite{OnTreesAndCharacters}, one can adapt the present proof for random $k$-cycles ($k$ fixed) at time $\left\lfloor \frac{1}{k}(n\log(n)+cn)\right\rfloor$, and we conjecture that the same limiting profile still holds for random conjugacy classes of size $o(n)$, as studied in \cite{BerestyckiSengul2014}, but that it would be technically much harder to adapt the present proof in that case. For this general case, a beautiful formula (Proposition 10.15 in \cite{LivreMelliot2017}) used in the proof of the Stanley-Féray formula, which allows to compute any reduced character as an expectation, $\chi^\lambda(\mu) = \bbE\cfg (-1)^{\text{inv}(\sigma_\mu)}\cfd$, might be very useful. \\
We conjecture that this profile also holds for the random involution walk studied by Megan Bernstein in  
\cite{Bernstein2018}, at time  $\left\lfloor\frac{1}{\log(p)}(\log(n)+c)\right\rfloor$.
For other problems where the limiting profile is known, see \cite{BayerDiaconis1992} and \cite{Lacoin2016}.
\subsection{Links with previous results and idea of the proof}\label{subsection12}
\paragraph{Links with previous results}
In 1981, Diaconis and Shahshahani showed in \cite{DiaconisShahshahani1981}, using representations of the symmetric group, a cutoff\footnote{In fact their lower bound is $1/e$ so it is not exactly a cutoff.} at $\left\lfloor\frac{1}{2}n\log(n)\right\rfloor$ for the random transposition shuffle, giving asymptotic inequalities at time $\left\lfloor\nlc\right\rfloor$, $c>0$ fixed. In 1987, Matthews, in \cite{Matthews1987}, refined these results thanks to a probabilistic proof. In 2011, Berestycki, Schramm and Zeitouni generalized in \cite{BerestyckiSchrammZeitouni2011} the previous result to the shuffle by random $k$-cycles, for $k$ fixed as $n \rightarrow \infty$, proving a cutoff at $\left\lfloor\frac{1}{k}n\log(n)\right\rfloor$, conjectured by Diaconis. Finally, in 2014, Berestycki and Şengül generalized again this result, in \cite{BerestyckiSengul2014}, to any conjugacy class whose support is $o(n)$, and without representation theory. \\
The proof in \cite{DiaconisShahshahani1981} relies on the so-called \textbf{Diaconis' upper bound lemma}, which leads to a sum over irreducible representations which they delicately bound with representation theory and analysis. Actually we can observe that the only place where a lot of information (we lose a factor $e$ in the limit $c\rightarrow\infty$ of the limit profile) is lost on the limit profile is at the very begining, when the Cauchy-Schwarz inequality is used in the proof of the upper bound lemma. Section \ref{section2} presents a remedy to this information loss, improving the upper bound lemma to an approximation lemma (Lemma \ref{approxlem}) which is asymptotically much more precise. Subsection \ref{subsection41}, quite technical, generalizes the asymptotic bounds of Diaconis and Shahshahani to any $c\in\bbR$. \\
Another crucial point of our proof is to pack together, in the sums over the irreducible representations $\lambda =(\lambda_1,...)$ of $\Sn$, all the partitions with the same $\lambda_1$. More precisely, Subsection \ref{subsection42} shows that when $j\in\bbN^*$ is fixed, we can study the sum over the partitions with $\lambda_1$ equal to $n-j$ as a sum over the partitions of the integer $j$, resulting in explicit manipulable formulas. \\
To understand where the limiting profile comes from, observe that, thanks to the lower bound of Matthews, the key observable is the number of fixed points. The limit profile is the distance between the asymptotic distribution of the number of fixed points of our walk at time $\left\lfloor\frac{1}{2}n\log(n)+cn\right\rfloor$, which is a $\text{Poiss}\left(1+e^{-2c}\right)$ distribution, and that of a pemutation taken uniformly at random, i.e. $\text{Poiss}(1)$. \\
\newline
Theorem \ref{theorem} stated above gives support to the following conjecture of Nathanaël Berestycki:

\begin{conjecture}\label{conj:CC}
  Let $\tau_n$ be the first time that all cards have been touched, and let $X_{\tau_n}$ be the state of the deck of cards at this (random) time. Then $\dvt (X_{\tau_n}, U_n) \to 0$ as $n \to \infty$.
\end{conjecture}
In other words, the conjecture says that $\tau_n$ is a stopping time at which the random permutation is well mixed for all practical purposes. Note that at time $\tau_n -1$ the permutation contains at least one fixed point, so that $\dvt (X_{\tau_n - 1},U_n)$ cannot converge to zero. Hence, the conjecture implies that $\tau_n$ is in some strong sense optimal for mixing the deck of cards. \\
Let us now explain in what way Theorem \ref{theorem} above is related to this conjecture. For any time $t$, let $G_t$ be the random graph which contains an edge $(i,j)$ if and only if the corresponding transposition has been applied at least once prior to time $t$. Then $G_t$ is essentially a
realisation of the Erd\H{o}s--Rényi 
random graph with parameters $n$ and $p = 1- \exp ( - t / {\binom{n}{2}})$. 
It is easy to check that any cycle of the random permutation $X_t$ at time $t$, considered as a set, is a subset of a connected component of $G_t$. 
Hence it makes sense to consider the cycle structure of the permutation \textit{restricted to any particular connected component of $G_t$}. Let $\mathfrak{C}_t$ be the largest component of $G_t$ (which is 
macroscopic if $t \ge cn$ for some $c>1$, and actually contains all vertices with high probability after time $\tau_n$). $\mathfrak{C}_t$ is called the giant component of $G_t$.
By a famous result of Schramm \cite{SchrammRT}, the distribution of the lengths of the largest cycles of $X_t$ within $\mathfrak{C}_t$, normalised by the total size $|\mathfrak{C}_t|$ of the giant component, converges to a Poisson--Dirichlet distribution (in the sense of finite dimensional distributions). Hence these largest cycles can be seen to coincide in the limit with the distribution of a uniform permutation on the giant component (see e.g. \cite{BerestyckiSchrammZeitouni2011}). A stronger version of Schramm's theorem would be the following conjecture (also by N. Berestycki):
\begin{conjecture}\label{conj:giant}
 Suppose $t \ge cn/2$ for some $c>1$. Given $\mathfrak{C}_t$, the distribution of $X_t|_{\mathfrak{C}_t}$, is approximately uniform, in the sense that
 $\dvt (X_t|_{\mathfrak{C}_t}, U_{\mathfrak{C}_t} ) \to 0$ in probability as $n \to \infty$, where $U_{\mathfrak{C}_t}$ is a uniform permutation on the giant component $\mathfrak{C}_t$.
\end{conjecture} 

 It is not hard to see that Conjecture \ref{conj:giant} implies Conjecture \ref{conj:CC}. Indeed, Conjecture \ref{conj:giant} implies a very precise description of the structure of $X_t$ close to the mixing time:
 if $t = \lfloor \tfrac12 n \log n + cn \rfloor$, then according to this conjecture $X_t$ would consist, if $\tau_n >t$ of a permutation that is approximately uniform on $n-1$ points, plus an extra fixed point; and would otherwise be indistinguishable from a uniform permutation if $\tau_n \ge t$. Such a description would imply that
 $$
 \dvt( X_t, U_n) = \dvt (\Fix (X_t), \Poiss(1)) + o(1),
 $$
 where $\Fix(X_t)$ is the number of fixed points of $X_t$.
 It is furthermore relatively easy to check that $\bbP( \tau_n > t ) \to e^{-2c}$ and hence, still assuming Conjecture \ref{conj:giant}, we would deduce
  $$
 \dvt( X_t, U_n) = \dvt (\Poiss(1+ e^{-2c}), \Poiss(1)),
 $$
 where the extra term $e^{-2c}$ in the right hand side accounts precisely for the probability that $\tau_n > t$. Of course, this last display is precisely the content of our Theorem \ref{theorem}.
\\

\paragraph{Organisation of the article}
In Section \ref{section2}, we present the improvement of Diaconis' upper bound lemma, using the non-commutative Fourier transform, which brings us back to group representations. In Section \ref{section3}, we will recall some results on the representations of the symmetric group, get precise estimations of the hook-length and Murnagham-Nakayama combinatorial formulas when the size $n$ of our partitions tend to infinity with $n-\lambda_1$ constant, and we will prove some some upper bounds useful in the sequel. In Section \ref{section4}, we will prove the announced theorem decomposing approximation by approximation. From now on, $k$ will denote without ambiguity the integer
\begin{equation*}
    k=k(n,c)= \left\lfloor\frac{1}{2}n\log(n)+cn\right\rfloor.
\end{equation*}

\paragraph{Idea of the proof} The algebraic objects $\widehat{\Sn}, \text{triv}, d_\lambda, s_\lambda$ and $\text{ch}^\lambda$ will be defined at the begining of Section \ref{section2}. For all $\sigma\in\Sn$, $\text{Fix}(\sigma)$ will denote the number of fixed points of the permutation $\sigma$. For $j\in\bbN^*$, let us also define the polynomial $T_j(z)$ by the formula $\summ{i=0}{j}\binom{z}{j-i}\frac{(-1)^i}{i!}$. The idea is to first fix $c\in\bbR$, and then to define for all $\epsilon>0$ an integer $M=M(c,\epsilon)$ such that when $n$ tends to infinity, all the following approximations are true up to $\epsilon$. \\
Rewriting the sum using the Fourier transform and the improvement of Diaconis' lemma,
\begin{equation*}
    \text{d}_{\text{1}}\left(P_n^{*k}, U_n\right) = \frac{1}{\abs{\Sn}}\sum_{\sigma\in \Sn}  \absc{\sum_{\lambda\in\widehat{\Sn}\backslash\ag \text{triv}\ad}d_\lambda s_\lambda^k\text{ch}^\lambda(\sigma)} \approx \frac{1}{n!}\sum_{\sigma\in \Sn} \absc{\sum_{\lambda_1\geq n-M} d_\lambda s_\lambda^k\text{ch}^\lambda(\sigma)}.
\end{equation*}
Then, thanks to the polynomial convergence lemma and letting $M\rightarrow\infty$, we will get
\begin{equation*}
    \frac{1}{n!}\sum_{\sigma\in \Sn} \absc{\sum_{\lambda_1\geq n-M} d_\lambda s_\lambda^k\text{ch}^\lambda(\sigma)} \approx \frac{1}{n!}\sum_{\sigma\in \Sn} \absc{\summ{j=1}{M}e^{-2jc}T_j(\text{Fix}(\sigma))} \approx \frac{1}{n!}\sum_{\sigma\in \Sn} \absc{\summ{j=1}{\infty}e^{-2jc}T_j(\text{Fix}(\sigma))}.
\end{equation*}
Finally, letting $n\rightarrow\infty$,
\begin{align*}
    \frac{1}{n!}\sum_{\sigma\in \Sn} \absc{\summ{j=1}{\infty}e^{-2jc}T_j(\text{Fix}(\sigma))}
    & = \frac{1}{n!}\sum_{\sigma\in \Sn} \absc{e^{-e^{-2c}}\left(1+e^{-2c}\right)^{\text{Fix}(\sigma)} - 1} \\
    & \approx \bbE \absc{e^{-e^{-2c}}\left(1+e^{-2c}\right)^{\text{Poiss}(1)} - 1} \\
    & = \text{d}_{1}\left(\text{Poiss}\left(1+e^{-2c}\right),\text{Poiss}(1)\right).
\end{align*}

\section{Improvement of Diaconis' upper bound lemma}\label{section2}
In this section we present the improvement of Diaconis' upper bound lemma. We will stay in the framework of finite groups, but this lemma can be used in a wider framework, of compact groups for example. Our aim is to get a better approximation than in \cite{DiaconisShahshahani1981} by not using Cauchy-Schwarz before Fourier. \\
Let $G$ be a finite group, $\bbC G$ the group algebra of $G$ and $\widehat{G}$ the set of the irreducible representations of $G$. We note $\text{triv}$ the trivial representation of $G$ and $\widehat{G}^*=\widehat{G}\backslash\ag \text{triv}\ad$. For $\alpha\in\widehat{G}$, we also name $\rho_\alpha$ the matrix of the representation $\alpha$, $\text{ch}^\alpha$ its character and $d_\alpha$ its dimension. 
Let us first recall the inversion formula for the non-commutative Fourier transform, well-explained in \cite{LivreMelliot2017}. For $f \du G \rightarrow \bbC$ and $g\in G$, we have 
\begin{equation*}
    f(g) = \sum_{\alpha\in\widehat{G}} \frac{d_\alpha}{\abs{G}} \text{Tr}(\rho^\alpha(g)^*\widehat{f}(\alpha)).
\end{equation*}
We deduce that for all $t\in\bbN$,
\begin{align*}
    \text{d}_{\text{1}}\left(P_n^{*t}, U_n\right)
    & = \sum_{g\in G} \abs{\Pet(g)-U(g)} \\
    & = \sum_{g\in G} \absc{\sum_{\alpha\in\widehat{G}}\frac{d_\alpha}{\abs{G}}\text{Tr}(\widehat{(\Pet-U)}(\alpha)\rho^\alpha(g)^*)} \\
    & = \sum_{g\in G} \absc{\sum_{\alpha\in\widehat{G}^*}\frac{d_\alpha}{\abs{G}}\text{Tr}(\widehat{\Pet}(\alpha)\rho^\alpha(g)^*)}.
\end{align*}
Besides, as $P$ is a function which is constant on every conjugacy class, we know that for each $\alpha$, by Schur's lemma, $\widehat{P}(\alpha)$ is a homothety, of ratio $s_\alpha = \frac{\text{Tr}(\widehat{P}(\alpha))}{d_\alpha}$. We hence obtain:
\begin{equation*}
    \text{d}_{\text{1}}\left(P_n^{*t}, U_n\right) = \frac{1}{\abs{G}}\sum_{g\in G} \absc{\sum_{\alpha\in\widehat{G}^*}d_\alpha s_\alpha^t\overline{\text{ch}^\alpha(g)}}.
\end{equation*}

Now, if instead of having a single group $G$ we have an increasing sequence of groups $(G_n)_\nn$, and if $t=t(n)$ is a well-chosen time depending on $n$ (and possibly on another parameter), we will wish to make $n$ tend to infinity inside our sums, and thus obtain a convergence to an explicit formula which will prove a cutoff or give a limiting profile. The idea of the following lemma is to spot a finite set of irreducible representations which will (asymptotically) have most of the mass, in order to approximate the sum over all irreducible representations by a sum over only finitely many terms, uniformly in $n$, and then be allowed to make $n$ tend to infinity inside the finite sum.
\begin{lemma}
\textbf{(Approximation lemma)}\label{approxlem} Let $G$ be a finite group and $S\subset \widehat{G}^*$.
Then:
\begin{equation*}
\absd{\text{d}_{\text{1}}\left(P_n^{*t}, U_n\right)-\frac{1}{\abs{G}}\sum_{g\in G}\absc{\sum_{\alpha\in S}d_\alpha s_\alpha^t\overline{\text{ch}^\alpha(g)}}} \leq \sum_{\alpha\in \widehat{G}^* \backslash S} d_\alpha \abs{s_\alpha}^t.
\end{equation*}
\end{lemma}
\paragraph{Proof}
Using the fact that $\absb{\abs{a}-\abs{b}}\leq \abs{a-b}$ and triangle inequalites,
\begin{align*}
    & \absd{\text{d}_{\text{1}}\left(P_n^{*t}, U_n\right)-\frac{1}{\abs{G}}\sum_{g\in G}\absc{\sum_{\alpha\in S}d_\alpha s_\alpha^t\overline{\text{ch}^\alpha(g)}}} \\
    \leq \; &  \frac{1}{\abs{G}}\sum_{g\in G}\absc{\sum_{\alpha\in \widehat{G} \backslash S}d_\alpha s_\alpha^t\overline{\text{ch}^\alpha(g)}} \\
    \leq \; & \frac{1}{\abs{G}}\sum_{g\in G}\sum_{\alpha\in \widehat{G} \backslash S}d_\alpha \abs{s_\alpha}^t\abs{\text{ch}^\alpha(g)} \\
    = \; &  \sum_{\alpha\in \widehat{G} \backslash S}d_\alpha \abs{s_\alpha}^t \frac{1}{\abs{G}}\sum_{g\in G}\abs{\text{ch}^\alpha(g)}.  \;\;\;\;\;\;\;\;\;\;\;\;\;\;\; (*)
\end{align*}
 Now, for every irreducible character $\alpha$, by Cauchy-Schwarz inequality and orthonormality of the characters,
 \begin{equation*}
     \frac{1}{\abs{G}}\sum_{g\in G}\abs{\text{ch}^\alpha(g)} \leq \frac{1}{\abs{G}} \sqrt{\abs{G}\sum_{g\in G}\abs{\text{ch}^\alpha(g)}^2} = 1.
 \end{equation*}
Plugging into $(*)$, this concludes the proof.

\section{The symmetric group and its representations}\label{section3}
\subsection{Hook-length formula}\label{subsection31}
We recall a few facts from the representation theory of the symmetric group, that we will naturally index by integer partitions $\lambda$. In a diagram associated to a partition, the hook of a box is the number of boxes which are above or on the right of our box. We call $\text{équ}(\lambda)$ the product of the hooks of the partition $\lambda$. For example, consider the partition $\lambda=(7,3,2,1,1)$ of the integer $14$ filled with its hooks: \\

\setlength{\unitlength}{0.75cm}
\begin{picture}(8,5)
\linethickness{0.2mm}
\put(0,0){\line(0,1){5}}
\put(1,0){\line(0,1){5}}
\put(2,0){\line(0,1){3}}
\put(3,0){\line(0,1){2}}
\put(4,0){\line(0,1){1}}
\put(5,0){\line(0,1){1}}
\put(6,0){\line(0,1){1}}
\put(7,0){\line(0,1){1}}

\put(0,0){\line(1,0){7}}
\put(0,1){\line(1,0){7}}
\put(0,2){\line(1,0){3}}
\put(0,3){\line(1,0){2}}
\put(0,4){\line(1,0){1}}
\put(0,5){\line(1,0){1}}

\put(6.35,0.3){1}
\put(5.35,0.3){2}
\put(4.35,0.3){3}
\put(3.35,0.3){4}
\put(2.35,0.3){6}
\put(1.35,0.3){8}
\put(0.23,0.3){11}

\put(2.35,1.3){1}
\put(1.35,1.3){3}
\put(0.35,1.3){6}

\put(1.35,2.3){1}
\put(0.35,2.3){4}

\put(0.35,3.3){2}

\put(0.35,4.3){1}

\put(8,0){.}

\end{picture} \\

In this case, we have:
\begin{equation*}
    \text{équ}(7,3,2,1,1)=11\times8\times6\times(4\times3\times2\times1)\times(6\times3\times1\times4\times1\times2\times1)=11\times8\times6\times4!\times\text{équ}(3,2,1,1).
\end{equation*}

We now recall the hook length formula, a proof of which can be found in Chapter 3 of \cite{LivreMelliot2017}.
\begin{prop}\label{prop1}\textbf{(Hook-length formula)}
If $\lambda$ is a partition of some integer $n$, then $d_\lambda = \frac{n!}{\text{équ}(\lambda)}$. In particular, $d_{(n-j,\lambda_2,\lambda_3,...)}\leq \binom{n}{j}d_{(\lambda_2,\lambda_3,...)}$.
\end{prop}

If $\lambda = (\lambda_1,\lambda_2,\lambda_3,...)$ is an integer partition, we will denote by $\lambda^*$ the truncated partition $(\lambda_2,\lambda_3,\lambda_4,...)$, where the largest row has been removed. For example if $\lambda = (n-7,3,2,1,1)$, $\lambda^*=(3,2,1,1)$ and in this case we have when $n \rightarrow \infty$,
\begin{equation*}
    d_\lambda = \frac{n!}{(n-7+4)(n-8+2)(n-9+1)(n-10)!}\frac{1}{\text{équ}(\lambda^*)}=\frac{n!}{(n-7)!\,\text{équ}(\lambda^*)}\left(1-\frac{7}{n}+O\left( \frac{1}{n^2}\right) \right).
\end{equation*}
This can be easily generalized and gives the following asymptotic formula:
\begin{prop}\label{prop2}\textbf{(Asymptotic hook-length formula)}
Let $j\geq 1$ and $\lambda_2, \lambda_3,...$ be fixed integers such that $\lambda_2+\lambda_3+... = j$. Then when $n\rightarrow\infty$,
\begin{equation*}
    d_{(n-j,\lambda_2,\lambda_3,...)}= \binom{n}{j}d_{(\lambda_2,\lambda_3,...)}\left( 1-\frac{j}{n}+O\left( \frac{1}{n^2}\right) \right).
\end{equation*}
\end{prop}
\paragraph{Proof}
Let $n\in \bbN^{*}$ and $\lambda = \lambda(n) = (n-j,\lambda_2,\lambda_3,...)$. Then when $n\rightarrow\infty$, denoting by $\lambda^{*'}$ the conjugated partition of the partition $\lambda^{*}=(\lambda_2,\lambda_3,...)$,
\begin{align*}
    d_{(n-j,\lambda_2,\lambda_3,...)} & = \frac{n!}{(n-j+\lambda^{*'}_1)(n-j-1+\lambda^{*'}_2)...(n-2j+1 +\lambda^{*'}_j)}\frac{1}{\text{équ}(\lambda_2,\lambda_3,...)} \\
    & = \frac{n!}{(n-j)!\,\text{équ}(\lambda_2,\lambda_3,...)}\frac{n-j}{n-j+\lambda^{*'}_1}\frac{n-j-1}{n-j-1+\lambda^{*'}_2}...\frac{n-2j+1}{n-2j+1+\lambda^{*'}_j} \\
    & =\frac{n!}{(n-j)!\,\text{équ}(\lambda_2,\lambda_3,...)}\left( 1-\frac{\lambda^{*'}_1}{n} +O\left( \frac{1}{n^2}\right) \right) ...\left( 1-\frac{\lambda^{*'}_j}{n}+O\left( \frac{1}{n^2}\right) \right)\\
    & = \binom{n}{j}d_{(\lambda_2,\lambda_3,...)}\left( 1-\frac{j}{n}+O\left( \frac{1}{n^2}\right) \right).
\end{align*}

\begin{remark}
Actually we will only need the equivalent, but the term in $-\frac{j}{n}$ allows us, in the next subsection, to have a better intuition of the modified character ratios.
\end{remark} 

\subsection{Character ratios}\label{subsection32}
Let $\tau$ be a transposition.
We define as in \cite{DiaconisShahshahani1981} the character ratio $r(\lambda)=\frac{\text{ch}^\lambda(\tau)}{d_\lambda}$. We can give different explicit formulas for this object, among which the following symmetric one, which follows from Lemma $7.14$ in \cite{LivreMelliot2017}. \\
If $\lambda = (\lambda_1,\lambda_2,...,\lambda_n)$ is a partition of the integer $n$, then we have:
\begin{equation*}
    r(\lambda)=\frac{1}{\binom{n}{2}}\summ{i=1}{n}\binom{\lambda_i}{2}-\binom{\lambda_i'}{2}.
\end{equation*}
The modified character ratio, as defined in Section \ref{section2}, writes as $s_\lambda = \frac{1}{n}+\frac{n-1}{n}r(\lambda)$ and takes into account that we pick the identity with probability $1/n$. The following upper bounds are given in \cite{LivreDiaconis1988}.
\begin{prop}\label{boundDiaconis}
If $\lambda$ is a partition of the integer $n$, then
\begin{equation*}
    s_\lambda \leq \frac{\lambda_1}{n} \,\,\;\;\text{ and }\,\,\;\;\abs{r(\lambda)} \leq \frac{\lambda_1}{n}.
\end{equation*}
 Moreover, if $\lambda_1 \geq \frac{n}{2}$, then
\begin{equation*}
    s_\lambda \leq 1- \frac{2(\lambda_1+1)(n-\lambda_1)}{n^2}.
\end{equation*}
\end{prop}

We will also need an asymptotic expansion of $s_\lambda$, easily obtainable from the explicit formula for $r(\lambda)$: If $j\in\bbN^*$ and $\lambda_2\geq\lambda_3\geq...\geq\lambda_r$ are non-negative integers such that $\lambda_2+...+\lambda_r=j$, then when $n\rightarrow\infty$, 
\begin{equation*}
    r(n-j,\lambda_2,...,\lambda_r) = \frac{1}{\binom{n}{2}}\left(\binom{n-j}{2}+O(1)\right) = 1 - \frac{2j}{n}+O\left(\frac{1}{n^2}\right),
\end{equation*}
and so 
\begin{equation*}
    s_\lambda = 1 - \frac{2j}{n}+O\left(\frac{1}{n^2}\right).
\end{equation*}
\begin{remark} In the general case, to guess a cutoff, we want to find a $t=t(n)$ for which $d_\alpha\abs{s_\alpha}^t=\theta(1)$ as $n\rightarrow \infty$, for the representations $\alpha$ which have the most mass. In the case of the symmetric group, as $d_\lambda \approx n^j$, we want to find $t$ such that $\abs{s_\lambda}^t \approx n^{-j}$. For instance, for random transpositions, it is very natural to expect a cutoff at $\frac{1}{2}n\log(n)$ from the formula of $s_\lambda$, as $\left(1-\frac{2j}{n}\right)^{\frac{1}{2}n\log(n)} \approx n^{-j}$. 
\end{remark}
\subsection{Mass transfer in the Young graph}\label{subsection33}
It will be convenient to use the formalism of the Young graph for some calculations. Here we are going to study, in the Young graph, a measure transfer from a row to the next one, which can be extended by recurrence to several lines. We will write $\lambda \vdash m$ for some $m\geq 1$ to indicate that $\lambda$ is a partition of the integer $m$. We will also write $\lambda \nearrow \Lambda$ if $\lambda \vdash m$ and $\Lambda \vdash m+1$ to say that the diagram of $\Lambda$ can be obtained from the diagram of $\lambda$ by adding a box. Let us fix an integer $j\geq 1$. We recall the transition formula for the dimensions of diagrams, which we can find in \cite{LivreKerov} or \cite{LivreMelliot2017}: if we fix $\lambda \vdash j$, then we have the following transfer, which may be of independent interest:
\begin{equation*}
    \sum_{\Lambda\du\lambda \nearrow \Lambda}d_\Lambda = (j+1) d_\lambda.
\end{equation*}
Let $j$ be an integer and $(\gamma_\lambda)_{\lambda\vdash j}$ a sequence of real numbers. We extend this line to the next line, $j+1$, as follows, following the edges of the graph: if $\Lambda\vdash j+1$, we set $\gamma_\Lambda=\sum_{\lambda\nearrow\Lambda}\gamma_\lambda$. Then we have the transfer:
\begin{prop}\label{masstransfer}
\begin{equation*}
    \sum_{\Lambda\vdash j+1}\gamma_\Lambda d_\Lambda = (j+1)\sum_{\lambda\vdash j}\gamma_\lambda d_\lambda.
\end{equation*}
\end{prop} 
\paragraph{Proof}
\begin{align*}
    \sum_{\Lambda \vdash j+1}\gamma_\Lambda d_\Lambda
    & = \sum_{\Lambda \vdash j+1} \left( \sum_{\lambda\du\lambda \nearrow \Lambda}\gamma_\lambda\right) d_\Lambda \\
    & = \sum_{\lambda\vdash j}\sum_{\Lambda\vdash j+1} \textbf{1}_{\lambda\nearrow\Lambda}\gamma_\lambda d_\Lambda \\
    & = \sum_{\lambda \vdash j} \gamma_\lambda \sum_{\Lambda\du\lambda \nearrow\Lambda}d_\Lambda \\
    & = (j+1) \sum_{\lambda\vdash j}\gamma_\lambda d_\lambda. \\
\end{align*}
\subsection{Permutations usually do not have only little cycles}\label{subsection34}
We set, for $n\in\bbN^*$ and $1\leq j \leq n$,
\begin{equation*}
    \mathfrak{S}_{n,j} = \ag \sigma \in\Sn \du \text{all the cycles of } \sigma \text{ are of length} \leq j\ad.
\end{equation*}
Let us show that when $j$ is fixed, $\mathfrak{S}_{n,j}$ is asymptotically much smaller than $\Sn$.

\begin{prop}\label{littlecycles}
Let $j\geq 2$ be a fixed integer. Then for $n$ large enough,
\begin{equation*}
    \log\left(\frac{\abs{\mathfrak{S}_{n,j}}}{\abs{\Sn}}\right)\leq\frac{-n\log(n)}{T(j)},
\end{equation*}
where $T(j) = 1+2+...+j$.
\end{prop} 
\paragraph{Proof}
We can see that in $\mathfrak{S}_{n,j}$, there are at most $(n+1)^j$ conjugacy classes, because such a conjugacy class is determined by the number of fixed points, $2$-cycles,..., $j$-cycles of a representative, each one necessarily between $0$ and $n$. Let us give an upper bound on the cardinality of such a class.
Let $n\geq j$ be a large integer, $\mu=(\mu_1,...,\mu_r)$ a partition of the integer $n$ such that $\mu_1\leq j$ and $\mu_r\geq 1$, and $\cC_\mu$ the associated conjugacy class.
Then if $k_q$ denotes the number of $\mu_i$ equal to $q$, we have for $n$ big enough: 
\begin{align*}
    \abs{\cC_\mu} & = \frac{n!}{2^{k_2}3^{k_3}...j^{k_j}k_2!k_3!...k_j!(n-2k_2-3k_3-...-jk_j)!} \\
    & \leq \frac{n!}{k_2!k_3!...k_j!(n-2k_2-3k_3-...-jk_j)!} \\
    & = \binom{n}{(2k_2,...,jk_j,(n-2k_2-...-jk_j))} \frac{(2k_2)!}{k_2!}...\frac{(jk_j)!}{k_j!}. \\
    & \leq \binom{n}{(\frac{n}{j},\frac{n}{j},...,\frac{n}{j})} \frac{(2k_2)!}{k_2!}...\frac{(jk_j)!}{k_j!}. \\
    & \leq j^n \frac{(2k_2)!}{k_2!}...\frac{(jk_j)!}{k_j!}.
\end{align*}
Moreover this latest product will be greater if the $k_i$ increase, so we can assume without loss of generality that $2k_2+...+jk_j \geq n-1$. One of the $k_i$ is therefore necessarily of cardinal greater than $\frac{n-1}{2+3+...+j}=\frac{n-1}{T(j)-1}$. Furthermore, as $(2k_2)!...(jk_j)!\leq n!$, we obtain:
\begin{equation*}
    \abs{\cC_\mu}\leq j^n\frac{n!}{\left(\frac{n-1}{T(j)-1}\right)!}.
\end{equation*}
Thus for $n$ large enough, 
\begin{equation*}
    \frac{\abs{\mathfrak{S}_{n,j}}}{\abs{\Sn}} \leq (n+1)^j j^n \frac{1}{\left(\frac{n-1}{T(j)-1}\right)!},
\end{equation*}
i.e.
\begin{equation*}
    \log\left(\frac{\abs{\mathfrak{S}_{n,j}}}{\abs{\Sn}}\right) \leq j\log(n+1) + n\log(j) -  \log\left(\left(\frac{n-1}{T(j)-1}\right)!\right) \sim - \frac{n\log(n)}{T(j)-1}.
\end{equation*}
As $T(j)-1<T(j)$, this leads to the desired asymptotic upper bound.

\begin{remark} This upper bound proves in particular that the ratio $\frac{\abs{\mathfrak{S}_{n,j}}}{\abs{\Sn}}$ tends to $0$, even multiplied by any power function, or polynomial. It is this fact that we will use. The case $j=1$ that we did not process is trivial because in this case $\abs{\mathfrak{S}_{n,1}}=1$. \\ 
Besides, if we had proceeded more carefully, we could have shown that $k_j \sim \frac{n}{j}$ maximizes the heavy terms of the cardinality of the conjugacy class, and therefore that $\log(\abs{\mathfrak{S}_{n,j}}) \sim \left(1-\frac{1}{j}\right)n\log(n) $.
\end{remark}
\subsection{Upper bound on the number of $q$-cycles}\label{subsection35}
For every permutation $\sigma\in\Sn$ and $q\in\bbN^*$, let $N_q(\sigma) = N_q^{(n)}(\sigma) $ denote the number of $q$-cycles in the cycle decomposition of $\sigma$. We recall the well-know law for the number of fixed points of a random permutation\footnote{For $m=0$, we apply the inclusion-exclusion principle to $\bigcup_{i=1}^n F_i$, where $F_i = \ag\sigma\in\Sn : \sigma(i)=i\ad$, and then generalize for any $m$.}
\begin{equation*}
    \bbP(\sigma\in\Sn \du N_1(\sigma)=m) = \frac{1}{m!}\summ{i=0}{n-m}\frac{(-1)^i}{i!} \; ,\; 0\leq m\leq n.
\end{equation*}
In particular, we deduce that for all $0\leq m \leq n$, $\bbP(\sigma\in\Sn \du N_1(\sigma)=m)\leq \frac{1}{m!}$.
Now we generalize this upper bound to the number of $q$-cycles.
\begin{prop}\label{majnombrecycles}
Let $q,m\in\bbN^*$, then
\begin{equation*}
    \bbP(\sigma\in\Sn \du N_q(\sigma) = m ) \leq \frac{1}{q^mm!}.
\end{equation*}
\end{prop}

\paragraph{Proof} 
As in the previous paragraph, if $\mu_i$ is a partition of the integer $n$, we denote by $k_q$ the number of $\mu_i$ equal to $q$.

\begin{align*}
& \bbP(\sigma\in\Sn \du N_q(\sigma) = m ) \\
= \; & \frac{1}{n!}\sum \limits_{\underset{k_q\,=\,m}{\mu\,\vdash\, n}} \abs{C_\mu} \\
= \; & \summ{r=1}{\infty}\sum \limits_{\underset{k_q\,=\,m}{\mu=(\mu_1,...,\mu_r\geq1)\,\vdash\, n}}\frac{1}{2^{k_2}3^{k_3}...r^{k_r}k_2!k_3!...k_r!(n-2k_2-3k_3-...-rk_r)!} \\
= \; & \frac{1}{q^mm!}\summ{r=1}{\infty}\sum \limits_{\underset{k_q\,=\,0}{\mu=(\mu_1,...,\mu_r\geq1)\,\vdash\, n-qm}}\frac{1}{2^{k_2}3^{k_3}...r^{k_r}k_2!k_3!...k_r!(n-2k_2-3k_3-...-rk_r)!} \\
\leq \; & \frac{1}{q^mm!}\summ{r=1}{\infty}\sum_{\mu=(\mu_1,...,\mu_r\geq1)\,\vdash\, n-qm}\frac{1}{2^{k_2}3^{k_3}...r^{k_r}k_2!k_3!...k_r!(n-2k_2-3k_3-...-rk_r)!} \\
= \; & \frac{1}{q^mm!}\bbP(\sigma\in\mathfrak{S}_{n-qm} \du N_q(\sigma) = 0 )\\
\leq \; & \frac{1}{q^mm!}.
\end{align*}


\section{Proof of Theorem \ref{theorem}}\label{section4}
For this whole section, we fix $c\in\bbR$. We recall that $k=k(n,c)=\left\lfloor\frac{1}{2}n\log(n)+cn\right\rfloor$.

\subsection{Bounding the error}\label{subsection41}
The upper bound is similar to the upper bound of the sum appearing in \cite{LivreDiaconis1988} after applying Diaconis' upper bound lemma. However, as we want a more precise result, there will be some additional technical difficulties as $c$ may be negative.

We can observe that the representations of the symmetric group which contribute the most in the sum
\begin{equation*}
    \text{d}_{\text{1}}\left(P_n^{*k}, U_n\right) = \frac{1}{\abs{\Sn}}\sum_{\sigma\in \Sn} \absc{\sum_{\lambda\in\widehat{\Sn}^*}d_\lambda s_\lambda^k\text{ch}^\lambda(\sigma)}
\end{equation*}
correspond to partitions with a large first row. We will therefore naturally split according to $\lambda_1$. We set for all $M\in\bbN^*$, and integer $n$ large enough,
\begin{equation*}
    S_M(n) = \ag\lambda\in\widehat{\Sn}^*\du \lambda_1\geq n-M\ad.
\end{equation*}
From Lemma \ref{approxlem}, we get that for all $M\geq1$,
\begin{equation*}
    \absd{\text{d}_{\text{1}}\left(P_n^{*k}, U_n\right)-\frac{1}{\abs{\Sn}}\sum_{\sigma\in \Sn}\absc{\sum_{\lambda\in S_M(n)}d_\lambda s_\lambda^k\text{ch}^\lambda(\sigma)}} \leq  \sum_{\lambda\in \widehat{\Sn} \;;\; \lambda_1 < n-M}  d_\lambda \abs{s_\lambda}^k.
\end{equation*}
It remains to prove that the right hand side of this inequality tends to $0$ uniformly in $n$ when $M\rightarrow\infty$, and to estimate the second term in the left hand side. Our first task is to bound the error in the approximation.

\begin{lemma}\label{majreste}\textbf{(Upper bound on the remainder)}\\
For all $\epsilon > 0$ there exist $M=M(c,\epsilon)\geq1$ and $n_0=n_0(M)\in\bbN$ such that if $n\geq n_0$, then
\begin{equation*}
     \sum_{\lambda_1 \leq n - M} d_\lambda |s_\lambda|^k \leq \epsilon.
\end{equation*}
\end{lemma}

\paragraph{Proof} We recall that $s_\lambda = \frac{1}{n}+\frac{n-1}{n}r(\lambda)$. Observe that if $\lambda$ is a partition of $n$ such that $r(\lambda)\geq 0$, then $r(\lambda') = - r(\lambda)$ and so $s_\lambda =\abs{s_\lambda} \geq \abs{s_{\lambda'}}$. Let us first bound $\sum_{\lambda_1 \leq n - 1}d_\lambda\abs{s_\lambda}^k$ splitting the sum into pieces. Note that $\lambda_1 = n$ corresponds to $r(\lambda) = 1$, i.e. to $\lambda = (n)$, the trivial representation, which disappeared when we used the Fourier transform. Likewise, $r(\lambda)=-1$ corresponds to $\lambda =(1^n)$.

\begin{align*}
    \sum_{r(\lambda) < 1} d_\lambda\abs{s_\lambda}^k 
    & = d_{(1^n)}\abs{s_{(1^n)}}^k+\sum_{-1<r(\lambda) \leq -\frac{2}{n}}d_\lambda\abs{s_\lambda}^k+\sum_{-\frac{2}{n}<r(\lambda) <\frac{2}{n}}d_\lambda\abs{s_\lambda}^k+\sum_{\frac{2}{n}\leq r(\lambda) < 1}d_\lambda\abs{s_\lambda}^k \\
    & = S_1 + S_2+S_3+S_4.
\end{align*}

Let us bound these different pieces separately. The first one is the easiest: 
\begin{equation*}
    S_1 = \left(1-\frac{2}{n}\right)^{\left\lfloor\nlc\right\rfloor} = o(1),
\end{equation*}
\begin{equation*}
    S_3 \leq \sum_{-\frac{2}{n}<r(\lambda) <\frac{2}{n}} d_\lambda\left(\frac{3}{n}\right)^{k} \leq \left(\sum_{\lambda\in \widehat{\Sn}^*} d_\lambda^2\right)\left(\frac{3}{n}\right)^{k} \leq n!\left(\frac{3}{n}\right)^{\nlc} = o(1),
\end{equation*}
\begin{equation*}
    S_2 = \sum_{\frac{2}{n}\leq r(\lambda) < 1}d_\lambda\left(|s_\lambda| -\frac{2}{n}\right)^k \leq \sum_{\frac{2}{n}\leq r(\lambda) < 1}d_\lambda |s_\lambda|^k \left(1-\frac{2}{n}\right)^{k} \leq e^{-\frac{2}{n}\left(\nlc\right)}S_4 = \frac{e^{-2c}}{n}S_4,
\end{equation*}
where we used in the upper bound for $S_2$ that $\abs{s_\lambda}\leq 1$. If we succeed in proving that $S_4$ is bounded (in $n$), then we will be able to conclude that $\sum_{r(\lambda) < 1} d_\lambda\abs{s_\lambda}^k$ is bounded (in $n$). We will bound a sum a little larger than $S_4$, namely $\sum_{0\leq r(\lambda) < 1}d_\lambda\abs{s_\lambda}^k$. Let us begin by a crude bound which will prove useful in the sequel. If $1\leq j \leq n$, we have 
\begin{equation*}
    \sum_{\lambda_1 = n-j}d_\lambda \leq \sum_{\lambda^*\vdash j}\binom{n}{j}d_{\lambda^*} \leq \binom{n}{j} \sqrt{(\sum_{\lambda^*\vdash j}1^2)(\sum_{\lambda^*\vdash j}d_{\lambda^*}^2)}\leq \frac{n^j}{j!}\sqrt{2^jj!}\leq \frac{n^j2^{j/2}}{\sqrt{j!}}, \;\;\;\;\;\;\; (**)
\end{equation*}
where the two first inequalities come from Proposition \ref{prop1} and Cauchy-Schwarz, and the before last inequality comes from the fact that each partitions of the integer $j$ can be seen as one of the $2^j$ subsets of the set with $j$ elements.
Therefore we have, using Proposition \ref{boundDiaconis} (note that $r(\lambda)\geq 0$ implies that $s(\lambda)>0$) 
\begin{align*}
    S_4 & \leq \sum_{{0\leq r(\lambda) < 1}} d_\lambda s_\lambda^k \\
    & =\summ{j=1}{n-1} \sum\limits_{\underset{0\leq r(\lambda) < 1}{\lambda_1 = n-j}} d_\lambda s_\lambda^k \\
    & \leq \summ{j=1}{\left\lfloor n/1000\right\rfloor}\left(\sum_{\lambda_1 = n-j} d_\lambda\right) \left(1-\frac{2j(n-j+1)}{n^2}\right)^k + \summ{j=\left\lfloor n/1000\right\rfloor+1}{n-1}\left(\sum_{\lambda_1 = n-j} d_\lambda\right)\left(1-\frac{j}{n}\right)^k \\
    & = A_1 + A_2.
\end{align*}

Let us bound $A_1$. We have, using $(**)$ and $1+x\leq \exp{x}$, 
\begin{align*}
    A_1
    & \leq \summ{j=1}{\left\lfloor n/1000\right\rfloor}\frac{n^j2^{j/2}}{\sqrt{j!}} \left(1-\frac{2j(n-j+1)}{n^2}\right)^k \\
    & \leq \summ{j=1}{\left\lfloor n/1000\right\rfloor} \frac{2^{j/2}}{\sqrt{j!}} e^{j\log(n)}e^{{-\frac{2j(n-j+1)}{n^2}\left(\nlc\right)}} \\
    & = \summ{j=1}{\left\lfloor n/1000\right\rfloor} \frac{2^{j/2}}{\sqrt{j!}} e^{j\log(n)}e^{-j\left(1-\frac{j-1}{n}\right)\left(\log(n)+2c\right)} \\
    & = \summ{j=1}{\left\lfloor n/1000\right\rfloor} \frac{2^{j/2}}{\sqrt{j!}} e^{-2jc}e^{j(j-1)\frac{\log(n)+2c}{n}}.
\end{align*}

Let $a_j(n)$ be the summand in the right hand side, and note that
\begin{equation*}
    \frac{a_{j+1}(n)}{a_j(n)} = \frac{e^{\frac{\log(2)}{2}-2c}}{\sqrt{j+1}}e^{2j\frac{\log(n)+2c}{n}}.
\end{equation*}
As a function of $j$ when $n$ is fixed, this is decreasing until $j=\frac{n}{4(\log(n)+2c)}$ and then increasing. If the first and the last ratios are (strictly) less than $1$, then we will have a subgeometric sum, which will hence be bounded.
The last ratio, at $\frac{n}{1000}$, is equal to
\begin{equation*}
    \sqrt{1000}e^{\frac{\log(2)}{2}-2c+\frac{4c}{1000}} n^{\frac{2}{1000}-\frac{1}{2}} \xrightarrow[n\rightarrow\infty]{} 0.
\end{equation*}
For the first ratio, we need to be a little more careful. At $j=1$, we can have a ratio much larger than $1$, all the more when $c$ is little (i.e negative and far from $0$). So we will need to split once more and consider the sum starting at a suitably chosen $M$, depending on $c$ but not on $n$. Thus, though the convergence is fast in the case of a positive $c$, already treated by Diaconis and Shahshahani, if $c$ is very negative, we will have to consider a very large amount of terms, and the convergence will be much slower.
Let $M$ be such that 
\begin{equation*}
    \frac{e^{\frac{\log(2)}{2}-2c}}{\sqrt{M+1}} \leq \frac{1}{4},
\end{equation*}
and $n$ large enough such that
\begin{equation*}
    e^{2M\frac{\log(n)+2c}{n}}\leq 2,
\end{equation*}
and that the ratio $\frac{a_{j+1}(n)}{a_j(n)}$ at $j=n/1000$ be less than $1/2$. Then as all the ratios from $j=M$ are less than $1/2$, we have: 
\begin{equation*}
    \summ{j=1}{n/1000}a_j(n) \leq \summ{j=1}{M}a_j(n) +a_M(n)\summ{i=1}{\infty}\frac{1}{2^i} \xrightarrow[n\rightarrow\infty]{} \summ{j=1}{M}\frac{2^{j/2}}{\sqrt{j!}} e^{-2jc} + \frac{2^{M/2}}{\sqrt{M!}} e^{-2Mc}.
\end{equation*}
Thus, as $c\in\bbR$ is fixed, $A_1$ is bounded uniformly in $n$.
Let us now treat $A_2$, which will be slightly easier.

We observe that for all $j\geq 0$, $j^j\leq j!3^j$, hence by $(**)$,

 \begin{equation*}
     \sum_{\lambda_1 = n-j}d_\lambda \leq\frac{n^j6^{j/2}}{j^{j/2}}.
 \end{equation*}
 
 Let $j$ be an integer between $n/1000$ and $n-1$. Then
 \begin{align*}
     \frac{n^j6^{j/2}}{j^{j/2}} \left(1-\frac{j}{n}\right)^k
     = \; & \frac{n^j6^{j/2}}{j^{j/2}} e^{k\log\left(1-\frac{j}{n}\right)} \\
     \leq \; & \frac{n^j6^{j/2}}{j^{j/2}} e^{-k\left(\frac{j}{n}+\frac{j^2}{2n^2}\right)} \\
     \leq \; & \frac{n^j6^{j/2}}{j^{j/2}} e^{-(\frac{1}{2}\log(n)+c)\left(j+\frac{n}{2\cdot 10^6}\right)} \\
    = \; & 6^{j/2} e^{\frac{j}{2}\log\left(\frac{n}{j}\right)} e^{-c\left(j+\frac{n}{2\cdot 10^6}\right)}e^{-\frac{1}{4\cdot 10^6}n\log(n)}\\
    \leq \; & 6^{j/2} e^{\frac{j}{2}\log(1000)} e^{\abs{c}\left(j+\frac{n}{2\cdot 10^6}\right)}e^{-\frac{1}{4\cdot 10^6}n\log(n)}\\
    \leq \; & e^{n\frac{\log(6)}{2}} e^{n\log(1000)} e^{\abs{c}\left(n+\frac{n}{2\cdot 10^6}\right)}e^{-\frac{1}{4\cdot 10^6}n\log(n)}\\
    = \; & e^{Kn-K'n\log(n)},
 \end{align*}
 where $K$ is a real constant and $K'$ is a positive constant.
 Thus,
 \begin{equation*}
     A_2 \leq n e^{Kn-K'n\log(n)} \xrightarrow[n\rightarrow\infty]{} 0.
 \end{equation*}

Now we are able to conclude, using the bounds in the proof for $A_1$. 
Let $\epsilon > 0$, and let $M = M(c,\epsilon)\geq 1$ such that  $\frac{e^{\frac{\log(2)}{2}-2c}}{\sqrt{M+1}} \leq \frac{1}{4}$ and $2\frac{2^{M/2}}{\sqrt{M!}}e^{-2Mc}<\epsilon$. Then for $n$ large enough,
\begin{align*}
        \sum_{\lambda_1\leq n-M} d_\lambda \abs{s_\lambda}^k 
        & \leq S_1+ S_2 + S_3 + \summ{j=M}{n/1000}a_j(n) + A_2 \\
        & \leq \summ{j=M}{n/1000}a_j(n) +o(1) \\
        & \leq a_M(n)\summ{i=0}{\infty}\frac{1}{2^i} + o(1) \\
        & \leq 2\frac{2^{M/2}}{\sqrt{M!}}e^{-2Mc} + o(1) \\
        & < \epsilon + o(1) \; \text{as } n\rightarrow\infty.
\end{align*}

\subsection{Polynomial convergence lemma}\label{subsection42}
We now start to estimate the main term.
\begin{lemma}\label{lemTj}
Let $\ell\in\bbN^*$. Then when $n\rightarrow\infty$,
\begin{align*}
\frac{1}{n!}\sum_{\sigma\in\Sn}\absc{\summ{j=1}{\ell}\sum_{\lambda_1=n-j}d_\lambda s_\lambda^k\text{ch}^\lambda(\sigma)} =\frac{1}{n!}\sum_{\sigma\in\Sn}\absc{\summ{j=1}{\ell}e^{-2jc}T_j(\text{Fix}(\sigma))} + o(1),    
\end{align*}
where we recall that
\begin{equation*}
    T_j(z) = \summ{i=0}{j}\binom{z}{j-i}\frac{(-1)^i}{i!}.
\end{equation*}
\end{lemma}
Let us first show how the polynomials $T_j$, a key element of the proof, arise naturally.
\begin{lemma}\label{lemTj2}
Let $j\in\bbN^*$ be a fixed integer, and $\sigma\in\Sn$ a permutation with at least one cycle of length greater\footnote{It still works for $\sigma \in \Sn\backslash\mathfrak{S}_{n,j-1}$.} than $j$ (i.e. $\sigma \in \Sn\backslash\mathfrak{S}_{n,j}$). Then 
\begin{equation*}
    \frac{1}{j!}\sum_{\lambda\in\widehat{\Sn}\du \lambda_1 = n-j}d_{\lambda^*}\text{ch}^\lambda(\sigma) = T_j(\text{Fix}(\sigma)).
\end{equation*}
\end{lemma}

\paragraph{Proof of Lemma \ref{lemTj2}}
This proof is combinatorial and strongly relies on the Murnagham-Nakayama rule.
We first consider $\sigma\in\Sn\backslash\mathfrak{S}_{n,j}$ as an indeterminate in $\text{ch}^\lambda(\sigma)$ and recall that, for any permutation $\sigma$ and $q\in\bbN^*$, $N_q(\sigma)$ is the number of $q$-cycles in the cycle decomposition of $\sigma$. 
For example, if $\lambda = (n-4,1,1,1,1)$ and $\sigma$ has a cycle of length greater than $4$, we have, using the Murnagham-Nakayama formula and writing $N_i$ for $N_i(\sigma)$, 
\begin{equation*}
    \ch^\lambda(\sigma) = \binom{N_1}{4} + N_3N_1 + \binom{N_2}{2}-N_4 -\left(\binom{N_1}{3} - N_2N_1+N_3\right)+\left(\binom{N_1}{2}-N_2\right) -N_1 +1.
\end{equation*}We can observe that $\text{ch}^\lambda(\sigma)$ is a polynomial in $N_1(\sigma)=\text{Fix}(\sigma),N_2(\sigma),..., N_j(\sigma)$. The key observation is that we will be able to compute everything when we take the sum at $\lambda_1 = j$ constant, and that our polynomial, which seemingly has $j$ indeterminates, will in reality be a polynomial in only one variable, $N_1(\sigma)$, the number of fixed points of $\sigma$. This comes from the orthogonality of some characters and the mass transfer (Proposition \ref{masstransfer}), which will make all the other terms cancel.
Let us give a little more details.
\\
For the polynomial algebra $\bbC\cfg z_1, z_2,... \cfd$, we will not use the canonical basis generated by the $z_i^j$, but rather the one generated by the $\binom{z_i}{j}$, better suited here.

Let $\sigma\in\Sn\backslash\mathfrak{S}_{n,j}$. If $\lambda$ is a partition of $n$ such that $\lambda_1=n-j$, then the coefficient of $\binom{N_1(\sigma)}{j}$ in $\text{ch}^\lambda(\sigma)$ is naturally the number 
of ways we can fill the Young diagram of $\lambda^*$ with all the numbers from $1$ to $j$ with line and column growth, i.e. the number of standard tableaux of $\lambda^*$, which is $d_{\lambda^*}=\text{ch}^{\lambda^*}(Id)$. \\
More generally, if $j_1,...,j_r\in\bbN$ are such that $j_1+2j_2+...+rj_r=j$, then the coefficient of \begin{equation*}
    \binom{N_1(\sigma)}{j_1}\binom{N_2(\sigma)}{j_2}...\binom{N_r(\sigma)}{j_r}
\end{equation*}
in $\text{ch}^\lambda(\sigma)$ is \begin{equation*}
    \text{ch}^{\lambda^*}(r^{j_r},...,2^{j_2},1^{j_1}).
\end{equation*}
Thus, by orthogonality of the characters, the coefficient of $\binom{N_1(\sigma)}{j_1}\binom{N_2(\sigma)}{j_2}...\binom{N_r(\sigma)}{j_r}$ in the sum \begin{equation*}
    \sum_{\lambda\in\widehat{\Sn}\du \lambda_1 = n-j} d_{\lambda^*} \text{ch}^\lambda(\sigma)
\end{equation*}
is
\begin{equation*}
    \sum_{\lambda\in\widehat{\Sn}\du \lambda_1 = n-j} d_{\lambda^*} \text{ch}^{\lambda^*}\left(r^{j_r},...,2^{j_2},1^{j_1}\right) = \sum_{\lambda\in\widehat{\Sn}\du \lambda_1 = n-j} \text{ch}^{\lambda^*}(Id) \text{ch}^{\lambda^*}\left(r^{j_r},...,2^{j_2},1^{j_1}\right) = 0.
\end{equation*}
By mass transfer, we can also observe that for $1\leq j'\leq j_1$, if $\sigma$ has at least $j'$ fixed points (if it has less, the coefficient is zero), the coefficient of 

\begin{equation*}
    \binom{N_1(\sigma)}{j_1-j'}\binom{N_2(\sigma)}{j_2}...\binom{N_r(\sigma)}{j_r}
\end{equation*} 
in the sum

\begin{equation*}
    \sum_{\lambda\in\widehat{\Sn}\du \lambda_1 = n-j} d_{\lambda^*} \text{ch}^\lambda(\sigma)
\end{equation*}
is $(-1)^{j'}$ times $j(j-1)...(j-j'+1)$ the coefficient of 

\begin{equation*}
    \binom{N_2(\sigma)}{j_2}...\binom{N_r(\sigma)}{j_r}
\end{equation*} 
in the sum 
\begin{equation*}
    \sum_{\lambda\in\widehat{\mathfrak{S}_{n-j'}}\du \lambda_1 = n-j+j'} d_{\lambda^*} \text{ch}^\lambda(\sigma'),
\end{equation*}
where $\sigma'$ has $j'$ less fixed points than $\sigma$, but as many $i$-cycles for each $i\geq2$, coefficient which is zero except when $j_2=...=j_r=0$, where it is equal to $1$. To summarize, we have shown that
\begin{equation*}
    \frac{1}{j!}\sum_{\lambda\in\widehat{\Sn}\du \lambda_1 = n-j}d_{\lambda^*}\text{ch}^\lambda(\sigma) =  \binom{N_1(\sigma)}{j} -\binom{N_1(\sigma)}{j-1}+\frac{1}{2}\binom{N_1(\sigma)}{j-2}+...+ \frac{(-1)^j}{j!}=T_j(\text{Fix}(\sigma)).
\end{equation*}
\paragraph{Proof of Lemma \ref{lemTj}}
Using the fact that $\absb{\abs{a}-\abs{b}}\leq \abs{a-b}$ and the triangle inequality,
\begin{align*}
    & \abs{\frac{1}{n!}\sum_{\sigma\in\Sn}\absc{\summ{j=1}{\ell}\sum_{\lambda_1=n-j}d_\lambda s_\lambda^k\text{ch}^\lambda(\sigma)}-\frac{1}{n!}\sum_{\sigma\in\Sn}\absc{\summ{j=1}{\ell}e^{-2jc}T_j(\text{Fix}(\sigma))}} \\
    \leq \; & \frac{1}{n!}\sum_{\sigma\in\Sn}\absc{\summ{j=1}{\ell}\sum_{\lambda_1=n-j}d_\lambda s_\lambda^k\text{ch}^\lambda(\sigma)-\summ{j=1}{\ell}e^{-2jc}T_j(\text{Fix}(\sigma))} \\
    \leq \; & \frac{1}{n!}\sum_{\sigma\in\Sn}\summ{j=1}{\ell}\absc{\left(\sum_{\lambda_1=n-j}d_\lambda s_\lambda^k\text{ch}^\lambda(\sigma)\right)-e^{-2jc}T_j(\text{Fix}(\sigma))}.
\end{align*}
Let us now split the sum on $\Sn$ into two parts, along $\mathfrak{S}_{n,\ell}$ and $\Sn\backslash\mathfrak{S}_{n,\ell}$, and let us bound each of these two sums separately. We begin by the sum on $\mathfrak{S}_{n,\ell}$. As in our sum $0 \leq s_\lambda\leq 1$ and $\text{ch}^\lambda(\sigma)\leq d_\lambda$,
\begin{align*}
    & \frac{1}{n!}\sum_{\sigma\in\mathfrak{S}_{n,\ell}}\summ{j=1}{\ell}\absc{\sum_{\lambda_1=n-j}d_\lambda s_\lambda^k\text{ch}^\lambda(\sigma)-e^{-2jc}T_j(\text{Fix}(\sigma))} \\
    \leq \; & \frac{1}{n!}\sum_{\sigma\in\mathfrak{S}_{n,\ell}}\summ{j=1}{\ell}\sum_{\lambda_1=n-j}\left(d_\lambda s_\lambda^k\abs{\text{ch}^\lambda(\sigma)}+\abs{e^{-2jc}T_j(\text{Fix}(\sigma))}\right) \\
    \leq \; & \frac{1}{n!}\sum_{\sigma\in\mathfrak{S}_{n,\ell}}\summ{j=1}{\ell}\sum_{\lambda_1=n-j}\left(d_\lambda^2+e^{-2jc}(\ell+1)n^\ell\right) \\
    \leq \; & \frac{1}{n!}\sum_{\sigma\in\mathfrak{S}_{n,\ell}}\summ{j=1}{\ell}\sum_{\lambda_1=n-j}\left(\left(\binom{n}{j}d_{\lambda^*}\right)^2+e^{-2jc}(\ell+1)n^\ell\right) \\
    \leq \; &  K(\ell,c)n^{2\ell}\frac{\abs{\mathfrak{S}_{n,\ell}}}{\abs{\Sn}} \;\;\;\;\;\;\;\;\;\; \text{using } \binom{n}{j}d_{\lambda^*} \leq \frac{n^j}{j!}d_{\lambda^*}\leq n^j\leq n^\ell \\
    = \; & o(1),
\end{align*}
where $K(\ell,c)$ is a constant depending only on $l$ and $c$.
Let us treat the second sum, which we rewrite using Lemma \ref{lemTj2}:
\begin{align*}
    & \frac{1}{n!}\sum_{\sigma\in\Sn\backslash\mathfrak{S}_{n,\ell}}\summ{j=1}{\ell}\absc{\left(\sum_{\lambda_1=n-j}d_\lambda s_\lambda^k\text{ch}^\lambda(\sigma)\right)-e^{-2jc}T_j(\text{Fix}(\sigma))} \\ 
    = \; & \frac{1}{n!}\sum_{\sigma\in\Sn\backslash\mathfrak{S}_{n,\ell}}\summ{j=1}{\ell}\absc{\sum_{\lambda_1=n-j}\left(d_\lambda s_\lambda^k -e^{-2jc}\frac{d_{\lambda^*}}{j!}\right)\text{ch}^\lambda(\sigma)} \\
    \leq \; & \frac{1}{n!}\sum_{\sigma\in\Sn\backslash\mathfrak{S}_{n,\ell}}\summ{j=1}{\ell}\sum_{\lambda_1=n-j}\abs{d_\lambda s_\lambda^k -e^{-2jc}\frac{d_{\lambda^*}}{j!}}\abs{\text{ch}^\lambda(\sigma)}.
\end{align*}
Let us observe that
\begin{equation*}
    d_{(n-j,\lambda_2,...,\lambda_r)} s_{(n-j,\lambda_2,...,\lambda_r)}^k -e^{-2jc}\frac{d_{(\lambda_2,...,\lambda_r)}}{j!}= O\left(\frac{1}{n}\right)
\end{equation*}
for all $1\leq j\leq \ell$ and $\lambda_2\geq \lambda_3\geq...\geq\lambda_r\geq 1$ such that $\lambda_2+...+\lambda_r =j$. (Note that there are only a finite number of such terms.)
We split the right hand side according to whether $\max(N_1(\sigma),...,N_\ell(\sigma))$ is larger or smaller than $n^{\frac{1}{2\ell}}$. On the one hand,
\begin{align*}
    & \frac{1}{n!}\sum \limits_{\underset{\max(N_1(\sigma),...,N_\ell(\sigma))\,\leq\, n^{1/(2\ell)}}{\sigma\in\Sn\backslash\mathfrak{S}_{n,\ell}}}\summ{j=1}{\ell}\sum_{\lambda_1=n-j}\abs{d_\lambda s_\lambda^k -e^{-2jc}\frac{d_{\lambda^*}}{j!}}\abs{\text{ch}^\lambda(\sigma)} \\
    = \; & O\left(\frac{1}{n}\right) \frac{1}{n!}\sum \limits_{\underset{\max(N_1(\sigma),...,N_\ell(\sigma))\,\leq\, n^{1/(2\ell)}}{\sigma\in\Sn\backslash\mathfrak{S}_{n,\ell}}}\summ{j=1}{\ell}\sum_{\lambda_1=n-j}K(\ell,c)\max(N_1(\sigma),...,N_\ell(\sigma))^\ell \\
     = \; & O\left(\frac{1}{n}\right) \frac{1}{n!}\sum \limits_{\underset{\max(N_1(\sigma),...,N_\ell(\sigma))\,\leq\, n^{1/(2\ell)}}{\sigma\in\Sn\backslash\mathfrak{S}_{n,\ell}}}\summ{j=1}{\ell}\sum_{\lambda_1=n-j}O\left(n^\frac{1}{2}\right) \\
     = \; & O\left(n^{-\frac{1}{2}}\right).
\end{align*}
On the other hand, 
\begin{align*}
    & \frac{1}{n!}\sum \limits_{\underset{\max(N_1(\sigma),...,N_\ell(\sigma))\,>\, n^{\frac{1}{2\ell}}}{\sigma\in\Sn\backslash\mathfrak{S}_{n,\ell}}}\summ{j=1}{\ell}\sum_{\lambda_1=n-j}\abs{d_\lambda s_\lambda^k -e^{-2jc}\frac{d_{\lambda^*}}{j!}}\abs{\text{ch}^\lambda(\sigma)} \\
    = \; & \frac{1}{n!}\sum \limits_{\underset{\max(N_1(\sigma),...,N_\ell(\sigma))\,>\, n^{\frac{1}{2\ell}}}{\sigma\in\Sn\backslash\mathfrak{S}_{n,\ell}}}\summ{j=1}{\ell}\sum_{\lambda_1=n-j}O\left(\frac{1}{n}\right)K(\ell,c)\max(N_1(\sigma),...,N_\ell(\sigma))^\ell \\
    \leq \; & \bbP(\sigma\in\Sn\du \max(N_1(\sigma),...,N_\ell(\sigma))> n^{\frac{1}{2\ell}})O\left(\frac{1}{n}\right)O\left(n^\ell\right) \\
    \leq \; & \summ{i=1}{\ell}\bbP\left(\sigma\in\Sn\du N_i(\sigma)> n^{\frac{1}{2\ell}}\right)O\left(\frac{1}{n}\right)O\left(n^\ell\right) \\
    = \; & O\left(\frac{1}{\left(n^\frac{1}{2\ell}\right)!}\right)O\left(\frac{1}{n}\right)O\left(n^\ell\right)\;\;\;\;\;\;\;\;\;\; \text{from Proposition } \ref{majnombrecycles} \\
    = \; & o(1).
\end{align*}
\subsection{Neglecting polynomials of high degree}\label{subsection43}
\begin{lemma}\label{lemhautdegré}
Let $\epsilon>0$. There exist $M_0 = M_0(\epsilon,c)$ such that for all $M\geq M_0$ and $n\in\bbN^*$, 
\begin{equation*}
    \abs{\frac{1}{n!}\sum_{\sigma\in\Sn}\absc{\summ{j=1}{M}e^{-2jc}T_j(\text{Fix}(\sigma))}-\frac{1}{n!}\sum_{\sigma\in\Sn}\absc{\summ{j=1}{\infty}e^{-2jc}T_j(\text{Fix}(\sigma))}}\leq \epsilon.
\end{equation*}
\end{lemma}

\paragraph{Proof}
Let $M,n\in\bbN^*$. Then we have, using again $\absb{\abs{a}-\abs{b}}\leq \abs{a-b},$
\begin{align*}
    & \abs{\frac{1}{n!}\sum_{\sigma\in\Sn}\absc{\summ{j=1}{M}e^{-2jc}T_j(\text{Fix}(\sigma))}-\frac{1}{n!}\sum_{\sigma\in\Sn}\absc{\summ{j=1}{\infty}e^{-2jc}T_j(\text{Fix}(\sigma))}} \\
    \leq \; &  \frac{1}{n!}\sum_{\sigma\in\Sn}\absc{\summ{j=M+1}{\infty}e^{-2jc}T_j(\text{Fix}(\sigma))} \\
    = \; &  \summ{r=0}{\infty}\bbP(\sigma\in\Sn \du N_1(\sigma)=r)\absc{\summ{j=M+1}{\infty}e^{-2jc}T_j(r)} \\
    \leq \; & \summ{r=0}{\infty}\frac{1}{r!}\summ{j=M+1}{\infty}e^{-2jc}\abs{T_j(r)}\;\;\;\;\;\;\;\;\;\; \text{from Proposition } \ref{majnombrecycles} \text{ again.}
\end{align*}
Now we observe that if $r\geq j$,
\begin{equation*}
    \abs{T_j(r)}=\abs{\summ{i=0}{j}\binom{r}{j-i}\frac{(-1)^i}{i!}}\leq\summ{i=0}{j}\binom{r}{j-i}\abs{\frac{(-1)^i}{i!}}\leq \summ{i=0}{j}\binom{r}{j-i}\leq 2^r,
\end{equation*}
and if $r\leq j$,
\begin{equation*}
    \frac{1}{r!}\abs{T_j(r)}=\frac{1}{r!}\abs{\summ{i=j-r}{j}\binom{r}{j-i}\frac{(-1)^i}{i!}}\leq \frac{1}{r!(j-r)!}\summ{i=j-r}{j}\binom{r}{j-i}\leq \frac{1}{\left(\left(\frac{j}{2}\right)!\right)^2}2^r.
\end{equation*}
We therefore conclude that
\begin{align*}
    & \summ{r=0}{\infty}\frac{1}{r!}\summ{j=M+1}{\infty}e^{-2jc}\abs{T_j(r)} \\
    = \; & \summ{j=M+1}{\infty}e^{-2jc} \summ{r=0}{j}\frac{1}{r!}\abs{T_j(r)} \;+ \; \summ{j=M+1}{\infty} \summ{r=j+1}{\infty}\frac{1}{r!}e^{-2jc}\abs{T_j(r)} \\
    \leq \; & \summ{j=M+1}{\infty}\frac{e^{-2jc}}{\left(\left(\frac{j}{2}\right)!\right)^2}\summ{r=0}{j}2^r \; + \; \summ{j=M+1}{\infty}\summ{r=j+1}{\infty}\frac{1}{r!}e^{2r\abs{c}}2^r \\
    \leq \; & \summ{j=M+1}{\infty}\frac{e^{-2jc}}{\left(\left(\frac{j}{2}\right)!\right)^2}2^{j+1} \; + \; \summ{j=M+1}{\infty}\summ{r=j+1}{\infty}\frac{1}{r!}e^{2r\abs{c}}2^r \\
    = \; & o(1)
\end{align*}
when $M\rightarrow\infty$. \\
Before proving the last approximation, let us rewrite the infinite sum inside the absolute values. Let us define
\begin{equation*}
    f_c \du x \mapsto e^{-e^{-2c}}\left(1+e^{-2c}\right)^x - 1.
\end{equation*}
\begin{prop}\label{propréécriture}

\end{prop}
Let $N\in\bbN$. Then 
\begin{equation*}
    \summ{j=1}{\infty}e^{-2jc}T_j(N) = f_c(N).
\end{equation*}
\paragraph{Proof}
We just need to make a change of variables and swap the two sums:
\begin{align*}
    \summ{j=1}{\infty}e^{-2jc}T_j(N) & = \summ{j=1}{\infty}\summ{i=0}{j}e^{-2jc}\binom{N}{j-i}\frac{(-1)^i}{i!} \\
    & = \summ{j=1}{\infty}\summ{i=0}{j}e^{-2jc}\binom{N}{i}\frac{(-1)^{j-i}}{(j-i)!} \\
    & = \summ{j=0}{\infty}\summ{i=0}{j}e^{-2jc}\binom{N}{i}\frac{(-1)^{j-i}}{(j-i)!} \;\; -1 \\
    & = \summ{i=0}{\infty}\binom{N}{i}e^{-2ic}\summ{j=i}{\infty}e^{-2(j-i)c}\frac{(-1)^{j-i}}{(j-i)!} \;\; -1 \\
    & = \summ{i=0}{N}\binom{N}{i}e^{-2ic}e^{-e^{-2c}} \;\; -1 \\
    & = e^{-e^{-2c}}\left(1+e^{-2c}\right)^{N} - 1.
\end{align*}
\subsection{Conclusion of the proof}\label{subsection44}
\begin{lemma}\label{lemconvpoiss}
When $n\rightarrow\infty$, we have: 
\begin{equation*}
    \frac{1}{n!}\sum_{\sigma\in \Sn} \left|f_c\left(N_1^{(n)}(\sigma)\right)\right| \xrightarrow[n\rightarrow\infty]{} \bbE \left|{f_c\left(\text{Poiss}(1)\right)}\right|,
\end{equation*}
where $\text{Poiss}(1)$ denotes the Poisson law of parameter $1$.
\end{lemma}

\paragraph{Proof}
As factorials grow much faster than exponentials, and hence than $f_c$, we have as $n\rightarrow\infty$,
\begin{align*}
    & \absd{\bbE \left|f_c\left(\text{Poiss}(1)\right)\right|-\frac{1}{n!}\sum_{\sigma\in \Sn} \left|f_c\left(N_1^{(n)}(\sigma)\right)\right|} \\
    =\; & \absd{\summ{r=0}{\infty}\frac{e^{-1}}{r!} \left|f_c(r)\right|-\summ{r=0}{n}\frac{1}{r!}\left(\summ{i=0}{n-r}\frac{(-1)^i}{i!}\right) \left|f_c(r)\right|} \\
    = \; & \absd{\summ{r=0}{n}\frac{1}{r!}\left(\summ{i=n-r+1}{\infty}\frac{(-1)^i}{i!}\right) \left|f_c(r)\right|+\summ{r=n+1}{\infty}\frac{e^{-1}}{r!}\left|f_c(r)\right|} \\
    = \; & o(1).
\end{align*}\\
We are now ready to combine all our estimates.
\paragraph{Proof of Theorem \ref{theorem}}
Let $\epsilon > 0$ and $M, n_0$ such that for $n\geq n_0$, all the approximations be true up to $\epsilon$. Let $n\geq n_0$. \\
From Lemma \ref{approxlem} and Lemma \ref{majreste},
\begin{equation*}
    \absd{\text{d}_{\text{1}}\left(P_n^{*k}, U_n\right)\;-\;\frac{1}{n!}\sum_{\sigma\in \Sn} \absc{\sum_{\lambda_1\geq n-M} d_\lambda s_\lambda^k\text{ch}^\lambda(\sigma)}}\leq \epsilon.
\end{equation*}
From Lemma \ref{lemTj}, 
\begin{equation*}
    \absd{\frac{1}{n!}\sum_{\sigma\in \Sn} \absc{\sum_{\lambda_1\geq n-M} d_\lambda s_\lambda^k\text{ch}^\lambda(\sigma)}\;-\;\frac{1}{n!}\sum_{\sigma\in \Sn} \absc{\summ{j=1}{M}e^{-2jc}P_j(N_1(\sigma))}}\leq \epsilon.
\end{equation*}
From Lemma \ref{lemhautdegré}, 
\begin{equation*}
    \absd{\frac{1}{n!}\sum_{\sigma\in \Sn} \absc{\summ{j=1}{M}e^{-2jc}P_j(N_1(\sigma))}\;-\;\frac{1}{n!}\sum_{\sigma\in \Sn} \absc{\summ{j=1}{\infty}e^{-2jc}P_j(N_1(\sigma))}} \leq \epsilon.
\end{equation*}
From Lemma \ref{lemconvpoiss},  
\begin{equation*}
    \absd{\frac{1}{n!}\sum_{\sigma\in \Sn} \absc{\summ{j=1}{\infty}e^{-2jc}P_j(N_1(\sigma))}\;-\; \bbE \left|f_c(\text{Poiss}(1))\right|} \leq \epsilon.
\end{equation*}
Consequently, by triangle inequalities,
\begin{equation*}
    \absd{\text{d}_{\text{1}}\left(P_n^{*k}, U_n\right)\;-\; \bbE \left|f_c(\text{Poiss}(1))\right|}\leq 4\epsilon.
\end{equation*}
Thus, we proved that for all $c\in\bbR$,
\begin{equation*}
    \text{d}_{\text{1}}\left(P_n^{*k}, U_n\right) \xrightarrow[n\rightarrow\infty]{} \bbE\left|f_c(\text{Poiss}(1))\right|.
\end{equation*}
To conclude, let us rewrite this expectation into the natural form of the wording:

\begin{align*}
   & \;\bbE\left|f_c(\text{Poiss}(1))\right| \\
   = \;& \summ{r=0}{\infty}\frac{e^{-1}}{r!} \absc{e^{-e^{-2c}}\left(1+e^{-2c}\right)^r - 1} \\
   = \;& \summ{r=0}{\infty} \absd{\frac{\left(e^{1+e^{-2c}}\right)^{-1}}{r!}\left(1+e^{-2c}\right)^r - \frac{e^{-1}}{r!}1^r} \\
   = \;& \text{d}_{1}\left(\text{Poiss}\left(1+e^{-2c}\right),\text{Poiss}(1)\right).
\end{align*}
This concludes the proof of Theorem \ref{theorem}.

\paragraph{Acknowledgements}
I am very thankful to my former professor and advisor, Justin Salez, who introduced me to mixing times and then took great care of me during my master thesis. I would also like to thank Nathanaël Berestycki, for his hospitality when he invited me to the University of Vienna, and for numerous helpful suggestions.


\begin{thebibliography}{9}

\bibitem{BayerDiaconis1992}
Dave Bayer, Persi Diaconis. Trailing the Dovetail Shuffle to its Lair. \textit{Ann. Appl. Prob.}, 2(2):294-313, 1992.

\bibitem{BerestyckiSchrammZeitouni2011}
Nathanaël Berestycki, Oded Schramm, Ofer Zeitouni. Mixing times for random $k$-cycles and coalescence-fragmentation chains. \textit{Ann. Probab.},39(5):1815-1843, 2011.

\bibitem{BerestyckiSengul2014}
Nathanaël Berestycki, Bati Şengül, Cutoff for conjugacy-invariant random walks on the permutation group, \textit{Probab. Theor. Rel. Fields}, to appear.

\bibitem{Bernstein2018}
Megan Bernstein, A random walk on the symmetric group generated by random involutions. \textit{Electronic Journal of Probability}, 2018.

\bibitem{BernsteinNestoridi2018}
Megan Bernstein, Evita Nestoridi, Cutoff for random to random card shuffle, submitted


\bibitem{BianeCartes2002}
Philippe Biane. Combien de fois faut-il battre un jeu de cartes? 
\textit{Gaz. Math.} No. 91, 4-10, 2002.


\bibitem{LivreDiaconis1988}
Persi Diaconis. \textit{Group representations in probability and statistics.} Institute of Mathematical Statistics Lecture Notes - Monograph Series, 11. Institute of Mathematical Statistics, Hayward, CA, 1988.


\bibitem{DiaconisShahshahani1981}
Persi Diaconis, Mehrdad Shahshahani. Generating a random permutation with random transpositions. \textit{Z. Wahrsch. Verw. Gebiete}, 57(2):159-179, 1981.


\bibitem{FreslonCutoffQuantique2018}
Amaury Freslon. Cut-off phenomenon for random walks on free orthogonal quantum groups, to appear in \textit{Probab. Theory Relat. Fields}, 2018. https://doi.org/10.1007/s00440-018-0863-8


\bibitem{OnTreesAndCharacters}
Avital Frumkin, Gordon James, Yuval Roichman. On Trees and Characters. \textit{Journal of Algebraic Combinatorics} (2003) 17: 323. https://doi.org/10.1023/A:1025052922664

\bibitem{LivreKerov}
Sergei V. Kerov. Asymptotic Representation Theory of the Symmetric Group and its Application in Analysis

\bibitem{Lacoin2016}
Hubert Lacoin. The cutoff profile for the simple exclusion process on the circle, \textit{Ann. Probab.}, 44(5), 3399-3430, 2016.

\bibitem{LulovPak2002} 
Nathan Lulov and Igor Pak. Rapidly mixing random walks and bounds on characters of the symmetric group. \textit{J. Algebraic Combin.}, 16(2):151–163, 2002.

\bibitem{MarkovChainsAndMixingTimes}
David A. Levin, Yuval Peres, and Elizabeth L. Wilmer. \textit{Markov chains and mixing times}.
American Mathematical Society, Providence, RI, 2009. With a chapter by James G. Propp and
David B. Wilson.

\bibitem{Matthews1987}
Peter Matthews, A strong uniform time for random transpositions, \textit{J. Theoret. Probab.} 1 (1988), 411–423.

\bibitem{LivreMelliot2017}
Pierre-Loïc Mélliot. Representation theory of the symmetric group

\bibitem{Salezcours}
Justin Salez, \textit{Temps de mélange des chaînes de Markov, Notes de cours}, 2018, \url{https://www.lpsm.paris/pageperso/salez/mixing.html}

\bibitem{SchrammRT}
Oded Schramm, Compositions of random transpositions, \textit{Israel J. Math.} 147 221–243.

\bibitem{Serre1977}
Jean-Pierre Serre. Linear Representations of Finite Groups, volume 42 of \textit{Graduate
Texts in Mathematics.} Springer–Verlag, 1977.

\end{thebibliography}
\end{document}